\documentclass[11pt, twoside]{article}
\usepackage{latexsym}
\usepackage{amsmath}
\usepackage{amssymb}
\usepackage[all]{xy}
\usepackage{amsfonts}
\usepackage{verbatim}
\usepackage{amsthm}
\usepackage{mathrsfs}
\usepackage{epsfig}
\usepackage{xy}
\usepackage{array}
\usepackage{stmaryrd}
\usepackage{graphicx,color}
\usepackage{xcolor}
\usepackage{tikz}
\usetikzlibrary{arrows,calc}
\usepackage{etex}
\usepackage{mathdots}
\usepackage{float}
\usepackage{graphics}
\usepackage{pdflscape}

\usepackage{anysize,hyperref}
\input xypic
\xyoption{all}

\usepackage[perpage,symbol]{footmisc}
\usepackage{setspace}
\def\A{\mathcal{A}}

\def\C{\mathscr{C}}

\def\dr{\ar@{->}[r]}

\def\X{\mathscr{X}}

\def\add{\mbox{add}}
\def\End{\mbox{End}}
\def\Ext{\mbox{Ext}}
\def\Hom{\mbox{Hom}}

\begin{document}
\baselineskip=15pt
\title{\Large{\bf $n$-abelian quotient categories}}
\medskip
\author{Panyue Zhou and Bin Zhu\footnote{Corresponding author. Panyue Zhou was Supported by the Hunan Provincial Natural Science Foundation of China (Grant No. 2018JJ3205).
Bin Zhu was supported by the NSF of China (Grants No.11671221).}}

\date{}

\maketitle
\def\blue{\color{blue}}
\def\red{\color{red}}

\newtheorem{theorem}{Theorem}[section]
\newtheorem{lemma}[theorem]{Lemma}
\newtheorem{corollary}[theorem]{Corollary}
\newtheorem{proposition}[theorem]{Proposition}
\newtheorem{conjecture}{Conjecture}
\theoremstyle{definition}
\newtheorem{definition}[theorem]{Definition}
\newtheorem{question}[theorem]{Question}
\newtheorem{remark}[theorem]{Remark}
\newtheorem{remark*}[]{Remark}
\newtheorem{example}[theorem]{Example}
\newtheorem{example*}[]{Example}
\newtheorem{condition}[theorem]{Condition}
\newtheorem{condition*}[]{Condition}
\newtheorem{construction}[theorem]{Construction}
\newtheorem{construction*}[]{Construction}

\newtheorem{assumption}[theorem]{Assumption}
\newtheorem{assumption*}[]{Assumption}

\baselineskip=17pt
\parindent=0.5cm

\begin{abstract}
\baselineskip=16pt
Let $\C$ be an $(n+2)$-angulated category with shift functor $\Sigma$ and $\X$ be a cluster-tilting subcategory of $\C$. Then we show that the quotient category $\C/\X$ is an $n$-abelian category.
If $\C$ has a Serre functor, then $\C/\X$ is equivalent to an $n$-cluster tilting subcategory of an abelian category $\textrm{mod}(\Sigma^{-1}\X)$. Moreover, we also prove that $\textrm{mod}(\Sigma^{-1}\X)$ is
Gorenstein of Gorenstein dimension at most $n$. As an application, we generalize recent results of Jacobsen-J{\o}rgensen and Koenig-Zhu.\\[0.5cm]
\textbf{Key words:} $(n+2)$-angulated categories; $n$-abelian categories; cluster-tilting subcategories; $n$-Gorenstein categories.\\[0.2cm]
\textbf{ 2010 Mathematics Subject Classification:} 18E30; 18E10.
\end{abstract}

\pagestyle{myheadings}
\markboth{\rightline {\scriptsize   Panyue Zhou and Bin Zhu}}
         {\leftline{\scriptsize $n$-abelian quotient categories}}

\section{Introduction}
Cluster-tilting theory sprouted from the categorification of Fomin-Zelevinsky's cluster algebras. It is used to construct abelian categories from some
triangulated categories. By Buan-Marsh-Reiten \cite[Theorem 2.2]{bmr} in cluster categories, by Keller-Reiten
\cite[Proposition 2.1]{kr} in the $2$-Calabi-Yau case, then by Koenig-Zhu \cite[Theorem 3.3]{kz} and Iyama-Yoshino \cite[Corollary 6.5]{iy} in the general case, one can pass from triangulated categories to abelian categories by factoring out any cluster tilting subcategory.

Recall that the notion of cluster-tilting subcategory, which was introduced by Keller-Reiten in \cite{kr}, see also \cite{kz, iy}:
Let $\C$ be a triangulated category with shift functor $\Sigma$. A subcategory $\X$ of $\C$ is called a \emph{cluster-tilting subcategory} of $\C$, if
it is satisfies the following conditions:

(1) $\X$ is contravariantly finite and covariantly finite in $\C$.

(2) $X\in\X$ if and only if $\Hom_{\C}(X,\Sigma\X)=0$, i.e. $\Hom_{\C}(X,\Sigma M)=0,$ for any $M\in \X$.

(3) $X\in\X$ if and only if $\Hom_{\C}(\X,\Sigma X)=0$,  i.e. $\Hom_{\C}(M,\Sigma X)=0,$ for any $M\in \X$.

An object $X$ is called cluster-tilting, if $\add(X)$ is a cluster-tilting,  where $\add(X)$ is the subcategory of $\C$ consisting of direct summands of direct sum of finitely many copies of $X$.

In fact, Koenig and Zhu \cite[Lemma 3.2]{kz} show that $\X$ is cluster-tilting if and only if $\Hom_{\C}(\X,\Sigma\X)=0$, i.e. $\Hom_{\C}(X,\Sigma Y)=0, \forall X, Y\in \X$, and for any object $C\in\C$, there exists a triangle
$X_0\to X_1\to C\to \Sigma X_0$
where $X_0,X_1\in\X$.

We introduce the notion of cluster tilting subcategory in an $(n+2)$-angulated
category.

\begin{definition}
Let $\C$ be an $(n+2)$-angulated category with an autoequivalence $\Sigma$ and $\X$ an additive subcategory of $\C$. $\X$ is called \emph{cluster-tilting} if

(1) $\Hom_{\C}(\X,\Sigma\X)=0$.

(2)  For any object $C\in\C$, there exists an $(n+2)$-angle
$$X_0\xrightarrow{~~}X_1\xrightarrow{~~}\cdots\xrightarrow{~~}X_{n-1}\xrightarrow{~~}X_n\xrightarrow{~~}C\xrightarrow{~~}\Sigma X_0$$
where $X_0, X_1,\cdots,X_{n}\in\X$.

An object $X$ is called \emph{cluster-tilting }in the sense of Opperman-Thomas \cite[Definition 5.3]{ot} if $\add(X)$ is clsuter-tilting.
\end{definition}

Recently, Jacobsen-J{\o}rgensen showed that the following result, which was obtained in a special case in the first part of \cite[Theorem 5.6]{ot}.

\begin{theorem}\emph{\cite[Theorem 0.6]{jj}}\label{thm0}
Let $\C$ be  a $k$-linear Hom-finite $(n+2)$-angulated category with split idempotents, where $k$ is an algebraically closed field and $n$ is a positive integer.  Assume that $\C$ has a Serre functor $\mathbb{S}$, that is, an autoequivlence for which there are natural
equivalences $$D\emph{\Hom}_{\C}(Z,Y)\simeq\emph{\Hom}_{\C}(Y,\mathbb{S}Z),$$ where $D(-)=\emph{\Hom}_{k}(-,k)$ is the $k$-linear duality functor and $Y,Z\in \C$.
Assume that $X$ is a cluster tilting object in $\C$ with $\Gamma=\emph{\End}_{\C}(X)$. Then

\emph{(1)} The (additive) quotient category $\C/\emph{\add}(X)$ is an $n$-abelian category.

\emph{(2)} $\C/\emph{\add}(X)$ is an $n$-cluster-titling subcategory of \emph{mod}$\Gamma$.

\emph{(3)} The endomorphism algebra $\Gamma$ is an $n$-Gorenstein algebra, that is, each injective module has projective dimension $\leq n$, and
each projective module has injective dimension $\leq n$.
\end{theorem}

Our  main result is a generalisation of Theorem \ref{thm0}, and \cite[Theorem 3.3]{kz} which can be recovered
by setting $n=1$.

\begin{theorem}\label{thm00} (See Theorem \ref{thm2} and Theorem \ref{thm3} for details)
Let $\C$ be an $(n+2)$-angulated category with split idempotents and $\X$ be cluster-tilting subcategory of $\C$. Then

\emph{(1)} The (additive) quotient category $\C/\X$ is an $n$-abelian category.

\emph{(2)} $\C/\X$ is an $n$-cluster-tilting subcategory of \emph{mod}$(\Sigma^{-1}\X)$.

\emph{(3)} If $\C$ is a $k$-linear Hom-finite with a Serre functor, then the category \emph{mod}$(\Sigma^{-1}\X)$ is a Gorenstein category, that is, each injective object has projective dimension $\leq n$, and
each projective object has injective dimension $\leq n$.

Note that if $\X=\emph{\add}(X)$, Theorem \ref{thm00} is precisely Theorem \ref{thm0}, and if $n=1$, Theorem \ref{thm00} is precisely Theorem 3.3 and Theorem 4.3 in \emph{\cite{kz}}.
\end{theorem}

The paper is organised as follows: In Section 2, we review some elementary definitions
that we need to use, including $n$-abelian categories and (right) $(n+2)$-angulated categories.
In Section 3,  we  provide a general framework for passing from $(n+2)$-angulated
categories to $n$-abelian categories by factoring out certain subcategories with some specific properties and
prove our main result by proving that cluster tilting subcategories satisfy the properties.

\section{Preliminaries}
In this section, we recall the definitions of $n$-abelian categoires and (right) $(n+2)$-angulated categories.
\subsection{$n$-abelian categories}
Let $\A$ be an additive category and $f\colon A\rightarrow B$ a morphism in $\A$. A \emph{weak cokernel} of $f$ is a morphism
$g\colon B\rightarrow C$ such that for any $X\in\A$ the sequence of abelian groups
$$\A(C,X)\xrightarrow{~g^\ast~}\A(B,X)\xrightarrow{~f^\ast~}\A(A,X)$$
is exact. Equivalently, $g$ is a weak cokernel of $f$ if $gf=0$ and for each morphism
$h\colon B\rightarrow X$ such that $hf=0$ there exists a (not necessarily unique) morphism
$p\colon C\rightarrow X$ such that $h=pg$. These properties are subsumed in the following
commutative diagram
$$\xymatrix{A\ar[r]^{f}\ar[dr]_{0} & B\ar[r]^{g}\ar[d]^{h}&C\ar@{-->}[dl]^{p}  \\
& X &}$$
Clearly, a weak cokernel $g$ of $f$ is a cokernel of $f$ if and only if $g$ is an epimorphism.
The concept of a \emph{weak kernel} is defined dually.

\begin{definition}\cite[Definiton 2.2]{j} and \cite[Definition 2.4]{l2}
Let $\A$ be an additive category and $f_0\colon A_0\rightarrow A_1$  a morphism in
$\A$. An $n$-\emph{cokernel} of $f_0$ is a sequence
$$(f_1,f_2,\cdots,f_{n})\colon A_1\xrightarrow{~f_1~}A_2\xrightarrow{~f_2~}\cdots \xrightarrow
{~f_{n-1}~}A_n\xrightarrow
{~f_{n}~}A_{n+1}$$
such that the induced sequence of abelian groups
$$\xymatrix{0\xrightarrow{~~}\A(A_{n+1},B)\xrightarrow{~~} \A(A_{n},B)\xrightarrow{~~}\cdots \xrightarrow{~~}
\A(A_{1},B)\xrightarrow{~~}
\A(A_{0},B)}$$
is exact for each object $B\in\A$. That is, the morphism $f_i$ is a weak cokernel of $f_{i-1}$ for all $i=1,2,\cdots,n-1$ and $f_{n}$ is a cokernel of $f_{n-1}$. In this case, we say the sequence
$$A_0\xrightarrow{~f_0~}A_1\xrightarrow{~f_1~}A_2\xrightarrow{~f_2~}\cdots \xrightarrow
{~f_{n-1}~}A_n\xrightarrow
{~f_{n}~}A_{n+1} $$
is \emph{right $n$-exact}.

We can define \emph{$n$-kernel} and \emph{left $n$-exact} sequence dually. The sequence
$$A_0\xrightarrow{~f_0~}A_1\xrightarrow{~f_1~}A_2\xrightarrow{~f_2~}A_3\xrightarrow{~f_3~}\cdots\xrightarrow
{~f_{n-1}~}A_n\xrightarrow{~f_n~}A_{n+1}$$
is called \emph{$n$-exact} if it is
both right $n$-exact and left $n$-exact.
\end{definition}

\begin{definition}\cite[Definiton 3.1]{j}\label{def0}
Let $n$ be a positive integer. An \emph{$n$-abelian category} is an additive category $\A$
which satisfies the following axioms:
\begin{itemize}
\item[(A0)] The category $\A$ has split idempotents.

\item[(A1)] Every morphism in $\A$ has an $n$-kernel and an $n$-cokernel.

\item[(A2)] For every monomorphism $f_0\colon A_0\to A_1$ in $\A$ there exists an $n$-exact sequence:
$$A_0\xrightarrow{~f_0}A_1\xrightarrow{~f_1~}A_2
\xrightarrow{~f_2~}\cdots\xrightarrow{~f_3~}
A_{n-1}\xrightarrow{~f_{n-1}~}A_{n}
\xrightarrow{~f_n~}A_{n+1}.$$

\item[(A2)$^{\textrm{op}}$] For every epimorphism $g_n\colon B_n\to B_{n+1}$ in $\A$ there exists an $n$-exact sequence:
$$B_0\xrightarrow{~g_0}B_1\xrightarrow{~g_1~}B_2
\xrightarrow{~g_2~}\cdots\xrightarrow{~g_3~}
B_{n-1}\xrightarrow{~g_{n-1}~}B_{n}
\xrightarrow{~g_n~}B_{n+1}.$$
\end{itemize}
\end{definition}

\begin{remark} $1$-abelian categories are precisely abelian categories in the usual sense. It is easy
to see that abelian categories have split idempotents; thus, if $n=1$, then axiom (A0) in
Definition \ref{def0} is redundant.
\end{remark}

\begin{definition}\cite{j}
Let $\A$ be a $n$-abelian category.
$P\in\A$ is called \emph{projective} if for every epimorphism $g\colon A\to B$ in $\A$, the sequence
$$\Hom_{\A}(P,A)\xrightarrow{\Hom_{\A}(P,\ g)}\Hom_{\A}(P,B)\xrightarrow{~~}0$$ is exact.
We say that $\A$ is \emph{projectively generated} if for every objects $A\in\A$ there exists a projective object $P\in \A$ and an epimorphism $f\colon P \to A$.
\end{definition}

\subsection{(Right) $(n+2)$-angulated categories}

Let $\mathcal{A}$ be an additive category with an endofunctor $\Sigma:\mathcal{A}\rightarrow\mathcal{A}$. An $(n+2)$-$\Sigma$-$sequence$ in $\mathcal{A}$ is a sequence of morphisms
$$A_0\xrightarrow{f_0}A_1\xrightarrow{f_1}A_2\xrightarrow{f_2}\cdots\xrightarrow{f_{n-1}}A_n\xrightarrow{f_n}A_{n+1}\xrightarrow{f_{n+1}}\Sigma A_0.$$
Its {\em left rotation} is the $(n+2)$-$\Sigma$-sequence
$$A_1\xrightarrow{f_1}A_2\xrightarrow{f_2}A_3\xrightarrow{f_3}\cdots\xrightarrow{f_{n}}A_{n+1}\xrightarrow{f_{n+1}}\Sigma A_0\xrightarrow{(-1)^n\Sigma f_0}\Sigma A_1.$$
A \emph{morphism} of $(n+2)$-$\Sigma$-sequences is  a sequence of morphisms $\varphi=(\varphi_0,\varphi_1,\cdots,\varphi_{n+1})$ such that the following diagram commutes
$$\xymatrix{
A_0 \ar[r]^{f_0}\ar[d]^{\varphi_0} & A_1 \ar[r]^{f_1}\ar[d]^{\varphi_1} & A_2 \ar[r]^{f_2}\ar[d]^{\varphi_2} & \cdots \ar[r]^{f_{n}}& A_{n+1} \ar[r]^{f_{n+1}}\ar[d]^{\varphi_{n+1}} & \Sigma A_0 \ar[d]^{\Sigma \varphi_0}\\
B_0 \ar[r]^{g_0} & B_1 \ar[r]^{g_1} & B_2 \ar[r]^{g_2} & \cdots \ar[r]^{g_{n}}& B_{n+1} \ar[r]^{g_{n+1}}& \Sigma B_0
}$$
where each row is an $(n+2)$-$\Sigma$-sequence. It is an {\em isomorphism} if $\varphi_1, \varphi_2, \varphi_3, \cdots, \varphi_n$ are all isomorphisms in $\mathcal{A}$.
\medskip

We recall the notion of a right $(n+2)$-angulated category from \cite[Definition 2.1]{l2}.
Compare \cite[Definition 2.1]{l2}, the condition (RN1) is slightly different from that in \cite{l2}, we don't assume that
the class $\Theta$ is closed under direct summands.

\begin{definition}\label{d1}
A {\em right} $n$-\emph{angulated category} is a triple $(\mathcal{A}, \Sigma, \Theta)$, where $\mathcal{A}$ is an additive category, $\Sigma$ is an endofunctor of $\mathcal{A}$, and $\Theta$ is a class of $(n+2)$-$\Sigma$-sequences (whose elements are called right $(n+2)$-angles), which satisfies the following axioms:
\begin{itemize}
\item[(RN1)]
\begin{itemize}
\item[(a)] The class $\Theta$ is closed under isomorphisms and direct sums.

\item[(b)] For each object $A\in\mathcal{A}$ the trivial sequence
$$0\rightarrow A\xrightarrow{1_A}A\rightarrow 0\rightarrow\cdots\rightarrow 0\rightarrow 0$$
belongs to $\Theta$.

\item[(c)] Each morphism $f_0\colon A_0\rightarrow A_1$ in $\A$ can be extended to a right $(n+2)$-angle: $$A_0\xrightarrow{f_0}A_1\xrightarrow{f_1}A_2\xrightarrow{f_2}\cdots\xrightarrow{f_{n-1}}A_n\xrightarrow{f_n}A_{n+1}\xrightarrow{f_{n+1}}\Sigma A_0.$$
\end{itemize}
\item[(RN2)] If an $(n+2)$-$\Sigma$-sequence belongs to $\Theta$, then its left rotation belongs to $\Theta$.

\item[(RN3)] Each solid commutative diagram
$$\xymatrix{
A_0 \ar[r]^{f_0}\ar[d]^{\varphi_0} & A_1 \ar[r]^{f_1}\ar[d]^{\varphi_1} & A_2 \ar[r]^{f_2}\ar@{-->}[d]^{\varphi_2} & \cdots \ar[r]^{f_{n}}& A_{n+1} \ar[r]^{f_{n+1}}\ar@{-->}[d]^{\varphi_{n+1}} & \Sigma A_0 \ar[d]^{\Sigma \varphi_0}\\
B_0 \ar[r]^{g_0} & B_1 \ar[r]^{g_1} & B_2 \ar[r]^{g_2} & \cdots \ar[r]^{g_{n}}& B_{n+1} \ar[r]^{g_{n+1}}& \Sigma B_0
}$$ with rows in $\Theta$ can be completed to a morphism of  $(n+2)$-$\Sigma$-sequences.

\item[(RN4)] Given a commutative diagram
$$\xymatrix{
A_0\ar[r]^{f_0}\ar@{=}[d] & A_1 \ar[r]^{f_1}\ar[d]^{\varphi_1} & A_2 \ar[r]^{f_2} & \cdots\ar[r]^{f_{n-1}} & A_{n}\ar[r]^{f_{n}} & A_{n+1} \ar[r]^{f_{n+1}} & \Sigma A_0\ar@{=}[d] \\
A_0\ar[r]^{g_0} & B_1 \ar[r]^{g_2}\ar[d]^{h_1} & B_2\ar[r]^{g_2} & \cdots\ar[r]^{g_{n-1}} & B_{n}\ar[r]^{g_{n}} & B_{n+1} \ar[r]^{g_{n+1}} & \Sigma A_0\\
& C_2\ar[d]^{h_2} & & & & & \\
& \vdots\ar[d]^{h_{n-1}} & & & & & \\
& C_{n}\ar[d]^{h_{n}} & & & & & \\
& C_{n+1}\ar[d]^{h_{n+1}} & & & & & \\
& \Sigma A_1 & & & & & \\
}$$
whose top rows and second column belong to $\Theta$. Then there exist morphisms $\varphi_i\colon A_i\rightarrow B_i\ (i=2,3,\cdots,n+1)$, $\psi_j\colon B_j\rightarrow C_j\ (j=2,3,\cdots,n+1)$ and $\phi_k\colon A_k\rightarrow C_{k-1}\ (k=3,4,\cdots,n+1)$ with the following two properties:

(I) The sequence $(1_{A_1},\varphi_1, \varphi_2,\cdots,\varphi_{n+1})$ is a morphism of $(n+2)$-$\Sigma$-sequences.

(II) The $(n+2)$-$\Sigma$-sequence
$$A_2\xrightarrow{\left(
                    \begin{smallmatrix}
                      f_2 \\
                      \varphi_2 \\
                    \end{smallmatrix}
                  \right)} A_3\oplus B_2\xrightarrow{\left(
                             \begin{smallmatrix}
                               -f_3 & 0 \\
                               \varphi_3 & -g_2 \\
                               \phi_3 & \psi_2 \\
                             \end{smallmatrix}
                           \right)}
 A_4\oplus B_3\oplus C_2\xrightarrow{\left(
                                       \begin{smallmatrix}
                                         -f_4 & 0 & 0 \\
                                         -\varphi_4 & -g_3 & 0 \\
                                         \phi_4 & \psi_3 & h_2 \\
                                       \end{smallmatrix}
                                     \right)}A_5\oplus B_4\oplus C_3$$
$$\xrightarrow{\left(
                                       \begin{smallmatrix}
                                         -f_5 & 0 & 0 \\
                                         \varphi_5 & -g_4 & 0 \\
                                         \phi_5 & \psi_4 & h_3 \\
                                       \end{smallmatrix}
                                     \right)}\cdots\xrightarrow{\scriptsize\left(\begin{smallmatrix}
             -f_{n} & 0 & 0 \\
             (-1)^{n+1}\varphi_{n} & -g_{n-1} & 0 \\
             \phi_{n} & \psi_{n-1} & h_{n-2} \\
             \end{smallmatrix}
             \right)}A_{n+1}\oplus B_{n}\oplus C_{n-1}$$
$$\xrightarrow{\left(
                                                       \begin{smallmatrix}
                                                         (-1)^{n+1}\varphi_{n+1} &-g_{n} &0 \\
                                                          \phi_{n+1}& \psi_{n}& h_{n-1} \\
                                                       \end{smallmatrix}
                                                     \right)}B_{n+1}\oplus C_{n}\xrightarrow{(\psi_{n+1},\ h_{n})}C_{n+1}\xrightarrow{\Sigma f_1\circ h_n}\Sigma A_2 \hspace{10mm}$$
belongs to $\Theta$, and $h_n\circ\psi_n=\Sigma f_1\circ g_n$.
   \end{itemize}
\end{definition}
The notion of a \emph{left $n$-angulated category} is defined dually.
\vspace{1mm}

If $\Sigma$ is an equivalence and the class $\Theta$ is closed under direct summands, it is easy to see that the converse of an axiom (RN2) also holds, thus the right
$n$-angulated category $(\A,\Sigma, \Theta)$ is an $n$-angulated category in the sense of Geiss-Keller-Oppermann \cite[Definition 1.1]{gko} and in the sense of Bergh-Thaule \cite[Theorem 4.4]{bt}. If $(\A,\Sigma, \Theta)$ is a right $n$-angulated category, $(\A,\Omega, \Phi)$
is a left $n$-angulated category, $\Omega$ is a quasi-inverse of $\Sigma$ and $\Theta=\Phi$, then $(\A,\Sigma, \Theta)$ is an $n$-angulated category.

\section{$n$-abelian quotient categories}

Let $\C$ be an additive category and $\X$ be a subcategory of $\C$. A morphism $f\colon A\to B$ in $\C$ is called $\X$-monic if for any object $X\in\X$, we have that
$$\Hom_{\C}(B,X)\xrightarrow{\Hom_{\C}(f,\ X)}\Hom_{\C}(A,X)\xrightarrow{~~}0$$
is exact. For an $\X$-monic morphism $f\colon A\to B$, if $B\in \X$, then we call $f$ is a left $\X-$approximation of $A$. If any object $A\in \C$ has a left $\X-$approximation, then $\X$ is called covariantly finite in $\C$. We can define $\X$-epic morphisms, right $\X-$approximation of an object $A$, and contravariantly finite subcategories dually.

 We denote by $\C/\X$
the category whose objects are objects of $\C$ and whose morphisms are elements of
$\Hom_{\C}(A,B)/\X(A,B)$ for $A,B\in\C$, where $\X(A,B)$ the subgroup of $\Hom_{\C}(A,B)$ consisting of morphisms
which factor through an object in $\X$.
Such category is called the (additive) quotient category
of $\C$ by $\X$. For any morphism $f\colon A\to B$ in $\C$, we denote by $\overline{f}$ the image of $f$ under
the natural quotient functor $\C\to\C/\X$.

\begin{lemma}\label{lem0}
Let $\C$ be an $(n+2)$-angulated category and $\X$ an additive subcategory of $\C$ with an autoequivalence $\Sigma:\C\rightarrow \C$. If for any object $A\in\C$, there exists an $(n+2)$-angle
$$A\xrightarrow{~f_0~}X_1\xrightarrow{~f_1~}X_2\xrightarrow{~f_2~}\cdots\xrightarrow{~f_{n-1}~}X_n\xrightarrow{~f_n~}B\xrightarrow{~f_{n+1}~}\Sigma A$$
with $f_0$ is a left $\X$-approximation of $A$ and $X_1,X_2,\cdots,X_n\in\X$. Then the quotient
category $\C/\X$ is a right $(n+2)$-angulated category with the following endofunctor and right $(n+2)$-angles:
\begin{itemize}
\item[\emph{(1)}] For any object $A\in\C$, we take an $(n+2)$-angle
$$A\xrightarrow{~f_0~}X_1\xrightarrow{~f_1~}X_2\xrightarrow{~f_2~}\cdots\xrightarrow{~f_{n-1}~}X_n\xrightarrow{~f_n~}\mathbb{G}B\xrightarrow{~f_{n+1}~}\Sigma A$$
with $f_0$ is a left $\X$-approximation of $A$ and $X_1,X_2,\cdots,X_n\in\X$. Then $\mathbb{G}$ gives a well-defined endofunctor of $\C/\X$.
\item[\emph{(2)}] For any $(n+2)$-angle
 $$A_0\xrightarrow{~g_0~}A_1\xrightarrow{~g_1~}A_2\xrightarrow{~g_2~}A_3\xrightarrow{~g_3~}\cdots\xrightarrow{~g_{n-1~}}A_n\xrightarrow{~g_n~}A_{n+1}\xrightarrow{~g_{n+1}~}\Sigma A_0$$
with $g_0$ is an $\X$-monic, take the following commutative diagram of $(n+2)$-angles.
$$\xymatrix{
A_0 \ar[r]^{g_0}\ar@{=}[d]& A_1 \ar[r]^{g_1}\ar@{-->}[d]^{\varphi_1} & A_2 \ar[r]^{g_2}\ar@{-->}[d]^{\varphi_2}  & \cdots \ar[r]^{g_{n-1}}& A_n \ar[r]^{g_n}\ar@{-->}[d]^{\varphi_n}&A_{n+1}\ar[r]^{g_{n+1}}\ar@{-->}[d]^{\varphi_{n+1}} & \Sigma A_0 \ar@{=}[d]\\
A_0 \ar[r]^{f_0}&X_1 \ar[r]^{f_1} & X_2 \ar[r]^{f_2}  & \cdots \ar[r]^{f_{n-1}} & X_n \ar[r]^{f_n}&\mathbb{G}B\ar[r]^{f_{n+1}}& \Sigma A_0\\
}$$
Then we have a complex
$$A_0\xrightarrow{~\overline{g_0}~} A_1\xrightarrow{~\overline{g_1}~}A_2\xrightarrow
{~\overline{g_2}~}\cdots\xrightarrow{~\overline{g_{n-1}}~}A_{n}\xrightarrow{~\overline{g_{n}}~}A_{n+1}\xrightarrow{~\overline{\varphi_{n+1}}~}\mathbb{G}A_0.$$ We define right $(n+2)$-angles in $\C/\X$ as the complexes which are isomorphic to complexes obtained in this way.
\end{itemize}
\end{lemma}

\proof  Since the proof is similar to \cite[Theorem 3.7]{l1}, we
omit it. See also \cite[Remark 3.8]{l2}  \qed

\begin{theorem}\label{thm1}
Let $\C$ be an $(n+2)$-angulated category with split idempotents and $\X$ an additive subcategory of $\C$. Consider the following conditions:
\begin{itemize}
\item[\emph{(a)}] For any object, there exist two $(n+2)$-angles:
$$\Sigma^{-1}A\xrightarrow{~~}X_0\xrightarrow{~~}X_1\xrightarrow{~~}\cdots\xrightarrow{~~}X_{n-1}\xrightarrow{~~}X_n\xrightarrow{~f~}A$$
with $f$ is a right $\X$-approximation of $A$ and $X_0, X_1,\cdots,X_{n}\in\X$;
and
$$A\xrightarrow{~g~}X'_1\xrightarrow{~~}X'_2\xrightarrow{~~}\cdots\xrightarrow{~~}X'_n\xrightarrow{~~}X'_{n+1}\xrightarrow{~~}\Sigma A$$
with $g$ is a left $\X$-approximation of $A$ and $X'_1,X'_2,\cdots,X'_{n+1}\in\X$,

\item[\emph{(b)}] For any $(n+2)$-angle
 $$A_0\xrightarrow{f_0}A_1\xrightarrow{f_1}A_2\xrightarrow{f_2}A_3\xrightarrow{f_3}\cdots\xrightarrow{f_{n-1}}A_n\xrightarrow{f_n}A_{n+1}\xrightarrow{f_{n+1}}\Sigma A_0$$
in $\C$, $f_0$ is an $\X$-monic if $\overline{f_{n}}$ is an epimorphism in $\C/\X$ and
$f_{n}$ is an $\X$-epic if $\overline{f_{0}}$ is a monomorphism in $\C/\X$.
\end{itemize}
If the conditions (a) and (b) hold,  then $\C/\X$ is an $n$-abelian category.
\end{theorem}

\proof Since $\C$ has split idempotents, obviously, $\C/\X$ has split idempotents.
Thus (A0) is satisfied.

By Lemma \ref{lem0} and its dual, we obtain that
$(\C/\X,\mathbb{G})$ is a right $(n+2)$-angulated category and
$(\C/\X,\mathbb{H})$ is a left $(n+2)$-angulated category.
By the constructions of $\mathbb{G}$ and $\mathbb{H}$,
 we know that $\mathbb{G}=0=\mathbb{H}$.

For any morphism $\overline{f_0}\colon A_0\to A_1$, take an $(n+2)$-angle
 $$A_0\xrightarrow{~g_0~}X'_1\xrightarrow{~~}X'_2\xrightarrow{~~}\cdots\xrightarrow{~~}X'_n\xrightarrow{~~}X'_{n+1}=\mathbb{G}A_0\xrightarrow{~~}\Sigma A_0$$
with $g_0$ is a left $\X$-approximation of $A$ and $X'_1,X'_2,\cdots,X'_{n+1}\in\X$.
It is easy to see that $\binom{f_0}{g_0}\colon A_0\to A_1\oplus X'_1$ is an $\X$-monic.
So there exists an $(n+2)$-angle
$$A_0\xrightarrow{\binom{f_0}{g_0}}A_1\oplus X'_1\xrightarrow{f_1}A_2\xrightarrow{f_2}A_3\xrightarrow{f_3}\cdots\xrightarrow{f_{n-1}}A_n\xrightarrow{f_n}A_{n+1}\xrightarrow{f_{n+1}}\Sigma A_0.$$
Thus we have the following commutative diagram
$$\xymatrix{
A_0 \ar[r]^{\binom{f_0}{g_0}\qquad}\ar@{=}[d]& A_1\oplus X'_1 \ar[r]^{\quad f_1}\ar@{-->}[d]^{\varphi_1} & A_2 \ar[r]^{f_2}\ar@{-->}[d]^{\varphi_2}  & \cdots \ar[r]^{f_{n-1}}& A_n \ar[r]^{f_n}\ar@{-->}[d]^{\varphi_n}&A_{n+1}\ar[r]^{f_{n+1}}\ar@{-->}[d]^{\varphi_{n+1}} & \Sigma A_0 \ar@{=}[d]\\
A_0 \ar[r]^{g_0}&X'_1 \ar[r]^{g_1} & X'_2 \ar[r]^{g_2}  & \cdots \ar[r]^{g_{n-1}} & X'_n \ar[r]^{g_n}&X'_{n+1}\ar[r]^{g_{n+1}}& \Sigma A_0\\
}$$
of $(n+2)$-angles in $\C$.
It follows that $$A_0\xrightarrow{~\overline{f_0}~}A_1\xrightarrow{~\overline{f_1}~}A_2\xrightarrow{~\overline{f_2}~}\cdots\xrightarrow{~\overline{f_{n-2}}~}A_{n-1}
\xrightarrow{~\overline{f_{n-1}}~}A_{n}\xrightarrow{~\overline{f_{n}}~}A_{n+1}\xrightarrow{~~~}0$$
right $(n+2)$-angle in $\C/\X$.
This shows that $(\overline{f_1},\overline{f_2},\cdots, \overline{f_{n+1}})$ is an $n$-cokernel of $\overline{f_0}$.

So any morphism in $\C/\X$ has $n$-cokernel. Dually we can show that $\C/\X$ has $n$-kernel.  Thus (A1) is satisfied.

If $\overline{f_0}$ is a monomorphism in $\C/\X$,
there exists a right $(n+2)$-angle
$$A_0\xrightarrow{~\overline{f_0}~}A_1\xrightarrow{~\overline{f_1}~}A_2
\xrightarrow{~\overline{f_2}~}\cdots\xrightarrow{~\overline{f_{n-2}}~}
A_{n-1}\xrightarrow{~\overline{f_{n-1}}~}A_{n}\xrightarrow{~\overline{f_{n}}~}
A_{n+1}\xrightarrow{~~~}0$$
in $\C/\X$.
Without loss of generality, we may assume that it is induced by an $(n+2)$-angle
$$A_0\xrightarrow{~f_0~}A_1\xrightarrow{~f_1~}A_2\xrightarrow{~f_2~}
\cdots\xrightarrow{~f_{n-2}~}A_{n-1}\xrightarrow{~f_{n-1}~}A_{n}\xrightarrow{~f_{n}~}A_{n+1}\xrightarrow{~f_{n+1}~}\Sigma A_{1}$$
in $\C$. By assumption (b), we have that $f_{n}$ is an $\X$-epic.
Thus the above $(n+2)$-angle induces a left $(n+2)$-angle
$$0\xrightarrow{~~} A_0\xrightarrow{~\overline{f_0}~} A_1\xrightarrow{~\overline{f_1}~}A_2\xrightarrow
{~\overline{f_2}~}\cdots\xrightarrow{~\overline{f_{n-2}}~}A_{n-1}\xrightarrow{~\overline{f_{n-1}}~}A_{n}\xrightarrow{~\overline{f_{n}}~}A_{n+1}$$
in $\C/\X$.
Hence we have that
$$A_0\xrightarrow{~\overline{f_0}~}A_1\xrightarrow{~\overline{f_1}~}A_2\xrightarrow
{~\overline{f_2}~}\cdots\xrightarrow{~\overline{f_{n-2}}~}A_{n-1}\xrightarrow{~\overline{f_{n-1}}~}A_{n}\xrightarrow{~\overline{f_n}~}A_{n+1}$$
is an $n$-exact sequence in $\C/\X$.  Thus (A2) is satisfied. Dually we can show that (A2)$^{\textrm{op}}$ is satisfied.

Therefore $\C/\X$ is an $n$-abelian category. \qed
\begin{remark}
In Theorem \ref{thm1}, if $n=1$, we recover a result of Li \cite[Theorem 5.4]{li}.
\end{remark}

Now we prove that cluster tilting subcategories satisfy the conditions (a) and (b) in Theorem \ref{thm1}.  Thus we obtain the first part of our main conclusions (Theorem \ref{thm00}).

\begin{theorem}\label{thm2}
Let $\C$ be an $(n+2)$-angulated category with split idempotents and $\X$ a cluster-tilting subcategory of $\C$. Then $\C/\X$ is an $n$-abelian category.
\end{theorem}

\proof By the definition of cluster-tilting subcategroies, condition (a) holds in Theorem \ref{thm1}.

For any $(n+2)$-angle
 $$A_0\xrightarrow{f_0}A_1\xrightarrow{f_1}A_2\xrightarrow{f_2}\cdots\xrightarrow{f_{n-1}}A_n\xrightarrow{f_n}A_{n+1}\xrightarrow{f_{n+1}}\Sigma A_0$$
in $\C$, we have $f_{n+1}f_n=0$ implies that  $\overline{f_{n+1}}\circ\overline{f_n}=0$.

If $\overline{f_{n}}$ is an epimorphism in $\C/\X$, we have $\overline{f_{n+1}}=0$, namely, there exist morphisms
$u\colon A_{n+1}\to X_0 $ and $v\colon X_0\to \Sigma A_0$ where $X_0\in\X$ such that $f_{n+1}=vu$.

Now we will prove that $f_0$ is an $\X$-monic. For any morphism $a\colon A_0\to X$ with $X\in\X$,
we have $\Sigma a\circ v\in\Hom_{\C}(X_0,\Sigma X)=0$.
It follows that $\Sigma a \circ f_{n+1}=(\Sigma a\circ v)u=0$. So there exists a morphism
$\Sigma b\colon \Sigma A_1\to \Sigma X$ such that $\Sigma b\circ(-1)^n\Sigma f=\Sigma a$
and then $a_0=(-1)^nb f_0$. This shows that $f_0$ is an $\X$-monic.

 Dually we can show that if $\overline{f_{0}}$ is a monomorphism in $\C/\X$,
then $f_{n}$ is an $\X$-epic. Therefore condition (b) holds in Theorem \ref{thm1}.
By Theorem \ref{thm1}, we have that $\C/\X$ is an $n$-abelian category.  \qed

\begin{lemma}\label{lem1}
Let $\C$ be an $(n+2)$-angulated category with split idempotents, let $\X$ be a cluster-tilting subcategory of $\C$ and $\A$ be the $n$-abelian quotient category of $\C$ by $\X$.
Let $f_n\colon A_{n-1}\to A_n$ be a morphism in $\C$ which is a part of an $(n+2)$-angle
 $$\Sigma^{-1}A_{n+1}\xrightarrow{f_0}A_0\xrightarrow{f_1}A_1\xrightarrow{f_1}A_2\xrightarrow{f_2}
 \cdots\xrightarrow{f_{n-1}}A_{n-1}\xrightarrow{f_{n}}A_n\xrightarrow{f_{n+1}}A_{n+1}.$$
Then $\overline{f_{n}}$ is an epimorphism in $\A$ if and only if $\overline{f_{n+1}}=0$;
$\overline{f_1}$ is a monomorphism in $\A$ if and only if $\overline{f_0}$.

In particular, if $A_{n+1}$ is in $\X$ , then $\overline{f_n}$ is an epimorphism;
if $\Sigma^{-1}A_{n+1}\in\X$, then  $\overline{f_1}$ is a monomorphism.
\end{lemma}

\proof  We give the proof for the statement about epimorphisms, the proof of the case of monomorphisms
is obtained dually.  We first will deal with a special case:
Assume that $f_n\colon A_{n-1}\to A_n$ be a morphism in $\C$ which is a part of an $(n+2)$-angle
 $$\Sigma^{-1}X\xrightarrow{f_0}A_0\xrightarrow{f_1}A_1\xrightarrow{f_1}A_2\xrightarrow{f_2}
 \cdots\xrightarrow{f_{n-1}}A_{n-1}\xrightarrow{f_{n}}A_n\xrightarrow{f_{n+1}}X$$
where $X\in\X$.
We will show that $\overline{f_n}$ is an epimorphism in $\A$.
Let $\varphi\colon A_n\to B$ with $\overline{\varphi}\circ \overline{f_n}=0$.
Then there exists an object $X'\in\X$ and morphisms $p\colon A_{n-1}\to X'$ and $q\colon X'\to B$
such that $\varphi f_n=qp$.
By the definition of cluster-tilting, for an object $B\in\C$, there exists an
$(n+2)$-angle
 $$X_{-1}\xrightarrow{g_0}X_0\xrightarrow{g_1}X_1\xrightarrow{g_1}X_2\xrightarrow{g_2}
 \cdots\xrightarrow{g_{n-1}}X_{n-1}\xrightarrow{g_{n}}B\xrightarrow{g_{n+1}}\Sigma X_{-1}$$
where $X_{-1},X_0,\cdots,X_{n-1}\in\X$.
Since $\Hom_{\C}(\X,\Sigma\X)=0$, we obtain that $g_n$ is a right $\X$-approximation of $B$.
So there exists a morphism $s\colon X'\to X_{n-1}$ such that $g_ns=q$.
It follows that $\varphi f_n=gq=g_nsp$. Thus we have the following commutative diagram
$$\xymatrix@R=1.2cm{
\Sigma^{-1}X\ar[r]^{\quad f_0}\ar@{-->}[d]^{\varphi_{-1}=0}&A_0 \ar[r]^{f_1}\ar@{-->}[d]^{\varphi_0}& A_1 \ar[r]^{f_2}\ar@{-->}[d]^{\varphi_1}\ar@{-->}[dl]^{t_1} & A_2 \ar@{-->}[dl]^{t_2} \ar[r]^{f_3}\ar@{-->}[d]^{\varphi_2}  & \cdots \ar[r]^{f_{n-1}}& A_{n-1} \ar[r]^{f_n}\ar[d]^{sp}&A_{n}\ar@{-->}[dl]^{t_n} \ar[r]^{f_{n+1}}\ar[d]^{\varphi} & \ar@{-->}[dl]^{t} X \ar@{-->}[d]^{0}\\
X_{-1}\ar[r]^{g_0}&A_0 \ar[r]^{g_1}&X_1 \ar[r]^{g_2} & X_2 \ar[r]^{g_3}  & \cdots \ar[r]^{g_{n-1}} & X_{n-1} \ar[r]^{g_n}&B\ar[r]^{g_{n+1}}& \Sigma X_{-1}\\
}$$
of $(n+2)$-angles in $\C$.
Since $\Hom_{\C}(\X,\Sigma\X)=0$, we have $\varphi_{-1}=0$ and then $\varphi_0f_0=0$.
So there exists  a morphism $t_1\colon A_1\to X_0$ such that $t_1f_1=\varphi_0$.
It follows that $(\varphi_1-g_1t_1)f_1=\varphi_1f_1-g_1\varphi_0=0$. So there exists a morphism
$t_2\colon A_2\to X_1$ such that $\varphi_1-g_1t_1=t_2f_2$.
Continuing this process, there exists a morphism $t_n\colon A_n\to X_{n-1}$
such that $sp-g_{n-1}t_{n-1}=t_nf_n$.
It follows that $(\varphi-g_nt_n)f_n=g_nsp-g_nt_nf_n=g_ng_{n-1}f_n=0$.
So there exists  a morphism $t\colon X\to B$ such that $\varphi-g_nt_n=tf_{n+1}$.
Hence $\overline{\varphi}=0$. This shows that $\overline{f_n}$ is an epimorphism in $\A$.

Now we turn to the general case: Assume that
 $f_n\colon A_{n-1}\to A_n$ is a morphism in $\C$ which is a part of an $(n+2)$-angle
 $$\Sigma^{-1}A_{n+1}\xrightarrow{f_0}A_0\xrightarrow{f_1}A_1\xrightarrow{f_1}A_2\xrightarrow{f_2}
 \cdots\xrightarrow{f_{n-1}}A_{n-1}\xrightarrow{f_{n}}A_n\xrightarrow{f_{n+1}}A_{n+1}$$
such that $f_{n+1}$ factors through $X$ with $X\in\X$.
Then, by completing the rightmost hand square, we get a commutative diagram with rows being
$(n+2)$-angles:
$$\xymatrix{
\Sigma^{-1}X\ar[r]\ar@{-->}[d]&B_0 \ar[r]\ar@{-->}[d]& B_1 \ar[r]\ar@{-->}[d] & A_2 \ar[r]\ar@{-->}[d]  & \cdots \ar[r]& B_{n-1} \ar[r]^{g_{n}}\ar@{-->}[d]^{\varphi_n}&A_{n}\ar[r]^{g_{n+1}}\ar@{=}[d]& X\ar[d]\\
\Sigma^{-1}A_{n+1}\ar[r]^{\quad f_0}&A_0 \ar[r]^{f_1}&A_1 \ar[r]^{f_2} & A_2 \ar[r]^{f_3}  & \cdots \ar[r]^{f_{n-1}} & A_{n-1} \ar[r]^{f_n}&A_n\ar[r]^{f_{n+1}}& A_{n+1}\\
}$$
As shown above, we have that $\overline{g_n}$ is an epimorphism in $\A$ since $\overline{g_{n+1}}=0$.
This implies that $\overline{f_n}$ is also an
epimorphism in $\A$ since $\overline{f_n}\circ \overline{\varphi_n}=\overline{g_n}$.

The converse is easy: Suppose $\overline{f_n}$ is an epimorphism and
 $\overline{f_{n+1}}\circ \overline{f_n}$= 0, then $\overline{f_{n+1}}=0$.   \qed
\medskip

As an immediate consequence, we obtain the following result.
\begin{corollary}\label{cor1}
Let $\C$ be an $(n+2)$-angulated category with split idempotents, let $\X$ be a cluster-tilting subcategory of $\C$ and $\A$ be the $n$-abelian quotient category of $\C$ by $\X$.
Let
$$A_0\xrightarrow{f_0}A_1\xrightarrow{f_1}A_2\xrightarrow{f_2}\cdots\xrightarrow{f_{n-1}}A_n\xrightarrow{f_n}A_{n+1}\xrightarrow{f_{n+1}}\Sigma A_0$$
be $(n+2)$-angle in $\C$.
\begin{itemize}
\item If $\overline{f_{n+1}}=0$, then $$A_0\xrightarrow{~\overline{f_0}~}A_1\xrightarrow{~\overline{f_1}~}A_2
\xrightarrow{~\overline{f_2}~}\cdots\xrightarrow{~\overline{f_{n-2}}~}
A_{n-1}\xrightarrow{~\overline{f_{n-1}}~}A_{n}\xrightarrow{~\overline{f_{n}}~}
A_{n+1}$$ is a right $n$-exact sequence in $\A$.

\item If $\overline{\Sigma^{-1} f_{n+1}}=0$, then
$$ A_0\xrightarrow{~\overline{f_0}~} A_1\xrightarrow{~\overline{f_1}~}A_2\xrightarrow
{~\overline{f_2}~}\cdots\xrightarrow{~\overline{f_{n-2}}~}A_{n-1}\xrightarrow{~\overline{f_{n-1}}~}A_{n}\xrightarrow{~\overline{f_{n}}~}A_{n+1}$$
is a left  $n$-exact sequence in $\A$.

\item If $\overline{f_{n+1}}=0=\overline{\Sigma^{-1} f_{n+1}}$, then $$A_0\xrightarrow{~\overline{f_0}~}A_1\xrightarrow{~\overline{f_1}~}A_2\xrightarrow
{~\overline{f_2}~}\cdots\xrightarrow{~\overline{f_{n-2}}~}A_{n-1}\xrightarrow{~\overline{f_{n-1}}~}
A_{n}\xrightarrow{~\overline{f_n}~}A_{n+1}$$
is an $n$-exact sequence in $\A$.
\end{itemize}
\end{corollary}

\begin{lemma}\label{lem2}
Let $\C$ be an $(n+2)$-angulated category with split idempotents. Suppose that $\X$ is a cluster-tilting subcategory of $\C$ and $\A$ is the $n$-abelian quotient category of $\C$ by $\X$.
Then an object $M$ of $\A$ is a projective object if and only if $M\in \Sigma^{-1}\X$. Dually an object $N$ of $\A$ is an injective object if and only if $N\in\Sigma\X$.
\end{lemma}

\proof We prove the first statement only, the second one is obtained dually.

Given $X\in\X$. For any epimorphism $\overline{\alpha}\colon A\to B$ in $\A$, and any morphism
$\overline{\beta}\colon X[-1]\to B$, let
$$A_0\xrightarrow{f_0}A_1\xrightarrow{f_1}A_2\xrightarrow{f_2}\cdots
\xrightarrow{f_{n-2}}A_{n-1}\xrightarrow{f_{n-1}}A\xrightarrow{\alpha}B\xrightarrow{u}\Sigma A_0$$
be the $(n+2)$-angle into which $\alpha$ is embedded.
Since $\overline{\alpha}$ is an epimorphism, we have $\overline{u}=0$. That is to say, there exists
an object $X_0\in\X$ and morphisms $s\colon B\to X_0$ and $t\colon X_0\to A_0$ such that
$u=ts$. It follows that $s\beta=0$ since $\Hom_{\C}(\Sigma^{-1}X, X_0)=0$.
Thus we have $u\beta=ts\beta=0$. So there exists a morphism $\gamma\colon X[-1]\to A$ such that
$\beta=\alpha\gamma$. Hence $\overline{\beta}=\overline{\alpha}\circ\overline{\gamma}$.
This shows that $\Sigma^{-1}X$ is projective in $\A$.

Conversely, assume that $M$ is a projective object in $\A$.
By the definition of cluster-tilting, there exists an $(n+2)$-angle
$$\Sigma^{-1}X_0\xrightarrow{~~}\Sigma^{-1}X_1\xrightarrow{~~}\Sigma^{-1}X_2\xrightarrow{~~}\cdots
\xrightarrow{~~}\Sigma^{-1}X_{n-1}\xrightarrow{~~}\Sigma^{-1}X_{n}\xrightarrow{~\alpha~}M\xrightarrow{~u~}X_0$$
where $X_0,X_1,\cdots,X_n\in\X$.
By Lemma \ref{lem1}, $\overline{\alpha}$ is an epimorphism in $\A$ since $\overline{u}=0$.
Therefore $\overline{\alpha}\colon \Sigma^{-1}X_{n}\to M$ splits. Hence $M\in\Sigma^{-1}\X$.  \qed

Recall that an abelian category with enough projectives and injectives is called $n$-Gorenstein
if all projective objects of this category have  injective dimension $\leq n$, and all injective
objects have projective dimension $\leq n$. The maximum of the injective dimensions of projectives and the projective dimensions of injectives is called \emph{Gorenstein dimension} of the category.

\begin{theorem}\label{thm3}
Let $\C$ be a $k$-linear Hom-finite $(n+2)$-angulated category with split idempotents, $\X$ be a cluster-tilting subcategory of $\C$ and $\A$ the $n$-abelian quotient category of $\C$ by $\X$. Then:
\begin{itemize}
\item[\emph{(1)}] The category \emph{mod}$(\Sigma^{-1}\X)$ is an abelian category.
\item[\emph{(2)}] The category $\A$ is an $n$-cluster-tilting subcategory of \emph{mod}$(\Sigma^{-1}\X)$.
\item[\emph{(3)}] If $\C$ has a Serre functor $\mathbb{S}$,
then the category \emph{mod}$(\Sigma^{-1}\X)$ is Gorenstein of Gorenstein dimension at most $n$.
\end{itemize}
\end{theorem}

\proof (1) This statements follows from the fact that $\Sigma^{-1}\X$ has weak kernels \cite{au}.

(2) By the definition of cluster-tilting subcategories, for any object $A\in\A$, there exists an $(n+2)$-angle
$$\Sigma^{-1}X_0\xrightarrow{~~}\Sigma^{-1}X_1\xrightarrow{~~}\Sigma^{-1}X_2\xrightarrow{~~}\cdots
\xrightarrow{~~}\Sigma^{-1}X_{n-1}\xrightarrow{~~}\Sigma^{-1}X_{n}\xrightarrow{~\alpha~}A\xrightarrow{~u~}X_0$$
where $X_0,X_1,\cdots,X_n\in\X$. By Lemma \ref{lem1}, we obtain that $\overline{\alpha}\colon \Sigma^{-1}X_{n}\to A $ is an epimorphism in $\A$ since $\overline{u}=0$.
By Lemma \ref{lem2}, we have that $\Sigma^{-1}X_{n}$ is a projective object.
Thus $\A$ is projectively generated $n$-abelian category.
By \cite[Theorem 1.3]{k}, we have that $\A$ is an $n$-cluster-tilting subcategory of mod$(\Sigma^{-1}\X)$.

(3) By Lemma \ref{lem2} and (2), we have that $\Sigma^{-1}\X$ is a full subcategory of mod$(\Sigma^{-1}\X)$  consisting of projective objects. By Lemma \ref{lem2}, Lemma 2.1 in \cite{jj} and (2), we have that $\mathbb{S}\Sigma^{-1}\X$ is a full subcategory of mod$(\Sigma^{-1}\X)$  consisting of injective objects.
It is easy to see that $\mathbb{S}\Sigma^{-1}\X=\Sigma\X$ if $\C$ has a Serre functor $\mathbb{S}$.
Note that $\C$ has a Serre functor, the duality functor $D=\Hom(-,k)$ induces a duality between mod$(\Sigma^{-1}\X)$
and $(\Sigma^{-1}\X)\textrm{mod}^{\textrm{op}}$, in particular,
mod$(\Sigma^{-1}\X)$ has enough projectives and enough injectives.

For any injective object $\Sigma X\in \textrm{mod}(\Sigma^{-1}\X)$ with $X\in\X$, by the definition of
cluster-tilting, there exists an $(n+2)$-angle
$$X\xrightarrow{~~}\Sigma^{-1}X_0\xrightarrow{~~}\Sigma^{-1}X_1\xrightarrow{~~}\Sigma^{-1}X_2\xrightarrow{~~}\cdots
\xrightarrow{~~}\Sigma^{-1}X_{n-1}\xrightarrow{~~}\Sigma^{-1}X_{n}\xrightarrow{~~}\Sigma X\xrightarrow{~~}X_0$$
where $X_0,X_1,\cdots,X_{n}\in\X$.
By Corollary \ref{cor1}, we have that
$$\Sigma^{-1}X_0\xrightarrow{~~}\Sigma^{-1}X_1\xrightarrow{~~}\Sigma^{-1}X_2\xrightarrow{~~}\cdots
\xrightarrow{~~}\Sigma^{-1}X_{n-1}\xrightarrow{~~}\Sigma^{-1}X_{n}\xrightarrow{~~}\Sigma X$$
is an $n$-exact sequence in $\A$. Of course, it is also a complex.
We obtain that the following short exact sequences in mod$(\Sigma^{-1}\X)$.
$$0\xrightarrow{~~}\Sigma^{-1}X_0\xrightarrow{~~}\Sigma^{-1}X_1\xrightarrow{~~}B_2\xrightarrow{~~}0$$
$$0\xrightarrow{~~}B_2\xrightarrow{~~}\Sigma^{-1}X_2\xrightarrow{~~}B_3\xrightarrow{~~}0$$
$$\vdots$$
$$0\xrightarrow{~~}B_{n-1}\xrightarrow{~~}\Sigma^{-1}X_{n-1}\xrightarrow{~~}B_n\xrightarrow{~~}0$$
$$0\xrightarrow{~~}B_n\xrightarrow{~~}\Sigma^{-1}X\xrightarrow{~~}\Sigma X\xrightarrow{~~}0$$
Since $\Sigma^{-1}X_0$ and $\Sigma^{-1}X_1$ are projective in mod$(\Sigma^{-1}\X)$ and
$$0\xrightarrow{~~}\Sigma^{-1}X_0\xrightarrow{~~}\Sigma^{-1}X_1\xrightarrow{~~}B_2\xrightarrow{~~}0$$
is an exact sequence in mod$(\Sigma^{-1}\X)$, we have that
$B_2$ has projective dimension $\leq 1$.

Applying the functor $\Hom(?,-)$ to the exact sequence
$$0\xrightarrow{~~}B_2\xrightarrow{~~}\Sigma^{-1}X_2\xrightarrow{~~}B_3\xrightarrow{~~}0,$$
we have the following exact sequence:
$$0=\Ext^{2}(\Sigma^{-1}X_2,-)\xrightarrow{~~}\Ext^{2}(\Sigma^{-1}B_2,-)\xrightarrow{~\simeq~} \Ext^{3}(B_3,-)\xrightarrow{~~} \Ext^{3}(\Sigma^{-1}X_2,-)=0.$$
This shows that $B_3$ has projective dimension $\leq 2$.
Inductively, we have that $\Sigma X$ has projective dimension $\leq n$.
Similarly, we can show that any projective object in mod$(\Sigma^{-1}\X)$ has injective dimension $\leq n$.

Therefore mod$(\Sigma^{-1}\X)$ is Gorenstein of Gorenstein dimension at most $n$. \qed

Panyue Zhou\\
College of Mathematics, Hunan Institute of Science and Technology, Yueyang, Hunan, 414006, People's Republic of China.\\
E-mail: \textsf{panyuezhou@163.com}\\[0.3cm]
Bin Zhu\\
Department of Mathematical Sciences, Tsinghua University, Beijing, 100084, People's Republic of
China.\\
E-mail: \textsf{bzhu@math.tsinghua.edu.cn}

\end{document}